\documentclass[a4paper, 12pt]{article}
\usepackage{amsmath}
\usepackage{amsthm} 
\usepackage{mathtools} 
\usepackage{amssymb} 
\usepackage{amsxtra} 
\usepackage{amsfonts} 
\usepackage{txfonts} 
\usepackage{enumerate} 
\usepackage{bbm} 
\usepackage{dsfont} 
\usepackage{bm} 
\usepackage{thm-restate}

\usepackage[T1]{fontenc} 
\usepackage{textcomp} 
\usepackage[utf8]{inputenc}

\newtheoremstyle{break}
{\topsep}{\topsep}%
{\itshape}{}%
{\bfseries}{}%
{\newline}{}%
\usepackage{hyperref}

\usepackage{moresize} 
\usepackage{pdfpages} 
\definecolor{mypurple}{HTML}{7f00d4}
\definecolor{mygreen}{HTML}{71a300}
\hypersetup{
	colorlinks=true,
	linkcolor=mypurple,
	filecolor=purple,
	citecolor=mygreen,
	urlcolor=mypurple,
}

\usepackage[colorinlistoftodos]{todonotes}

\usepackage{marginnote}

\usepackage{graphicx} 
\usepackage{caption} %
\usepackage{subcaption}

\usepackage{pgf,tikz,pgfplots}

\usepackage{tikz-cd}

\usepackage{float} 
\usepackage{turnstile} 
	\usepackage{thmtools}

	\theoremstyle{break}
	\newtheorem{Teo}[subsection]{Theorem}
	\newtheorem*{mainthm}{Main Theorem}

	\theoremstyle{definition}
	\newtheorem{Def}[subsection]{Definition} 
	
	\theoremstyle{remark}
	\newtheorem{Remark}[subsection]{Remark}
	
	\theoremstyle{plain}
	\newtheorem{question}[subsection]{Question}
	\newtheorem{Lema}[subsection]{Lemma}
	\newtheorem{Prop}[subsection]{Proposition}
	\newtheorem{Coro}[subsection]{Corollary}

		\newenvironment{Dem}{\noindent \vspace{-0.1cm}{\it Proof:}}{\hfill $\square$ \vspace{0.5 cm}} 
	\newenvironment{Demof}[1]{\noindent \vspace{-0.1cm}{\it Proof {#1}:}}{\hfill $\square$ \vspace{0.5 cm}} 
	\newenvironment{Not}{\vspace{0.3cm}\noindent \textsc{Notation:}}{\vspace{0.3cm}}{}

	\newcommand{\concept}[1]{\textbf{#1}} 
	
	\newcommand{\forc}[2]{\dststile{#2}{#1}}
	
	\newcommand{\tup}[1]{\langle#1\rangle}

	\newcommand{\RR}{\mathbb{R}}
	\newcommand{\QQ}{\mathbb{Q}}
	\newcommand{\cc}{\mathfrak{c}}

	\DeclareMathOperator{\lh}{lh}
	\DeclareMathOperator{\dom}{dom}
	\DeclareMathOperator{\im}{im}
	
	\DeclareMathOperator{\field}{field}
	\DeclareMathOperator{\genpol}{gp}
	\DeclareMathOperator{\expfield}{exp}

	\usepackage{setspace} 
	
	\title{The transcendence degree of the reals over certain set-theoretical subfields}
	\author{Azul Fatalini\\ Ralf Schindler}
	\date{}
	
\begin{document}
		\maketitle

\begin{abstract}
	It is a well-known result that, after adding one Cohen real, the transcendence degree of the reals over the ground-model reals is continuum. We extend this result for a set $X$ of finitely many Cohen reals, by showing that, in the forcing extension, the transcendence degree of the reals over a combination of the reals in the extension given by each proper subset of $X$ is also maximal. This answers a question of Kanovei and Schindler \cite{Kanovei2020}. 	
\end{abstract}

{\def\thefootnote{}\footnotetext{The authors have been supported by the Deutsche Forschungsgemeinschaft (DFG, German Research Foundation) under the Excellence Strategy EXC 2044–390685587, Mathematics Münster: Dynamics–Geometry–Structure.}}

This article investigates a folklore result about the forcing extension by one Cohen real: the transcendence degree of the reals over the set of reals in the ground model is of cardinality $\mathfrak{c}$ (in the extension), where $\mathfrak{c}$ is the cardinality of the continuum.
	This is obtained by noticing that one Cohen real can be split in a perfect set of Cohen reals which are mutually generic and, therefore, algebraically independent over the ground model reals. 
	 
	It turns out it is possible to generalize this statement in two directions: taking into account other type of reals added (observed independently by Ben de Bondt and Elliot Glazer, see Corollary \ref{coro: trans degree for any inner model with different reals}) and adding more than one Cohen real (see Theorem \ref{th: transcendence degree many cohen}), which is the main result of this paper:\\
	
		\begin{mainthm}[Theorem \ref{th: transcendence degree many cohen}]
		Let $X$ be a finite set of mutually generic Cohen reals over $V$.  
		In $V[X]$, consider the minimum field $F\subseteq \RR$ such that $F\supseteq \bigcup_{Y\subsetneq X} \RR^{V[Y]}$. 
		Then, in $V[X]$ the transcendence degree of $\RR$ over $F$ is continuum. 
	\end{mainthm}
	\vspace{0.5cm}
	
	In \cite{Kanovei2020}, the authors asked whether a certain set-theoretical subfield of reals is a proper subfield in the
	Cohen model. 
	Theorem \ref{th: transcendence degree many cohen} shows that it is a proper subfield and, moreover, the reals have
	transcendence degree continuum over this subfield.
	
	Actually, one of the main obstacles to prove Theorem \ref{th: transcendence degree many cohen} is to show that the field $F$ is a proper subfield of $\RR$. 
	In the case of $n=1$, i.e., in the folklore result, this is for free: the Cohen real is a \emph{new real}, namely it does not belong to the ground model reals, witnessing that the corresponding field $F$ is a proper subfield.
	
	In the case of $n=2$, and $X=\{x,y\}$, it is tempting to take $x\oplus y$ (the real produced by alternating the digits of $x$ and $y$, considering them as elements in $\omega^\omega$) as a witness for $F$ being a proper subfield of $\RR$. 
	Nevertheless, after a closer examination, one can see that $x\oplus y$ actually belongs to the field $F$, since $x\oplus y = (x\oplus \vec{0})+(\vec{0}\oplus y)$. 
	The proof of Theorem \ref{th: transcendence degree many cohen} shows that an appropriate witness is the \emph{composition} of the reals in $X$ as functions from $\omega$ to $\omega$. \\
	
	This paper is organized as follows. 
	In Section \ref{section: Prereq}, we will go through the basic prerequisites of field theory and real analysis needed for the proofs of the auxiliary results that will lead to our main theorem.
	Section \ref{section: context} contains the proof of the first direction of generalization of the folklore result, when one is allowed to add any real, not necessarily Cohen. 
	Finally, Section \ref{section: the theorems} contains the proof of the main theorem, generalizing the folklore result in a second direction: allowing more than one Cohen real. 
	We conclude by including a result about the relationship between the transcendence bases of $\RR^{V[x,y]}$, $\RR^{V[x]}$ and $\RR^{V[y]}$ over $\RR^V$.

	\section{Prerequisites}\label{section: Prereq}
	\subsection{Field theory}\label{subsection: prereq algebra}
	
	\begin{Def}\label{def: prereq field extension, field generated}
		A field $F$ is an \concept{extension} of a field $K$ if $K$ is a subfield of $F$, namely, 
		$K\subseteq F$ and the operations on $K$ are the ones on $F$ restricted to $K$. 
		
		If $F$ is a field and $X\subseteq F$, then the \concept{subfield} (resp. \concept{subring}) \concept{generated by $X$} is the intersection of all subfields (resp. subrings) of $F$ that contain $X$. 
		If $F$ is an extension of $K$ and $X\subseteq F$, then the subfield (resp. subring) generated by $K\cup X$ is denoted by $K(X)$ (resp. $K[X]$). 
	\end{Def}
	
	\begin{Not}
		Given a field $K$ and $n\in \omega$, we denote by $K[x_0, \dots, x_{n-1}]$ the ring of polynomials in $n$ variables over $K$.
	\end{Not}
	
	First, let us recall the construction of a field extension for a given set $X$ of generators.
	
	\begin{Teo} [{\cite[Theorem 1.3 Chapter V]{Hungerford2003}}] \label{th: prereq description field generated by X}
		If $F$ is an extension field of a field $K$ and $X\subseteq F$, then the subfield $K(X)$ consists of all elements of the form 
		\[\frac{f(u_0, \dots, u_{n-1})}{g(u_0, \dots, u_{n-1})}= f(u_0, \dots, u_{n-1})g(u_0, \dots, u_{n-1})^{-1}, \]
		where $n\in \omega$, 
		$f,g\in K[x_0, \dots , x_{n-1}]$, $u_0, \dots, u_{n-1}\in X$ and $g(u_0, \dots, u_{n-1})\neq 0$.
	\end{Teo}
	
	\begin{Def} 

\label{def: prereq algebraic or trans element}
		Let $F$ be an extension field of $K$.  
		An element $u$ of $F$ is said to be \concept{algebraic} over $K$ iff $u$ is a root of some nonzero polynomial $f\in K[x]$. 
		If $u$ is not a root of any nonzero $f\in K[x]$, $u$ is \concept{transcendental} over $K$. 
	\end{Def}
	
	\begin{Def}\label{def: prereq alg closed and closure}
		A field $F$ is \concept{algebraically closed} if every non-constant polynomial $f\in F[x]$ has a root in $F$. 
		If $F$ is an extension field of $K$, $F$ is algebraically closed, and every element of $F$ is algebraic over $K$, we say that $F$ is an \concept{algebraic closure} of $K$.
	\end{Def}
	
	\begin{Remark}
		If $K$ is a subfield of $\RR$, we will often talk about \emph{the} algebraic closure of $K$ and it will mean the algebraic closure of $K$ that is a subfield of $\mathbb{C}$. 
		Moreover, we will be interested in the algebraic closure of $K$ relative to $\RR$, which is the algebraic closure of $K$ ($\subseteq \mathbb{C}$) intersected with $\RR$. 
		We will denote the algebraic closure of $K$ relative to $\RR$ simply by $\overline{K}$.
	\end{Remark}
	
	\begin{Def}\label{def: prereq alg independent}
		Let $F$ be an extension field of $K$ and $S$ a subset of $F$. 
		$S$ is \concept{algebraically dependent} over $K$ if for some $n<\omega$ there is a nonzero polynomial $f\in K[x_0, \dots, x_{n-1}]$ such that $f(s_0,\dots , s_{n-1})=0$ for some distinct $s_0, \dots, s_{n-1}\in S$. 
		$S$ is \concept{algebraically independent} over $K$ if $S$ is not algebraically dependent over $K$.
	\end{Def} 
	
	\begin{Def}\label{def: prereq trans base}
		Let $F$ be an extension field of $K$. 
		A \concept{transcendence base} of $F$ over $K$ is a subset $S$ of $F$ which is algebraically independent over $K$ and is a maximal (with respect to $\subseteq$) set with this property. 
	\end{Def}

		Let $F$ be an extension field of $K$. 
		Let $S$ be a transcendence base of $F$ over $K$.
		Then every transcendence base of $F$ over $K$ has the same cardinality as $S$ 
		\cite[Theorems 1.8-1.9 Chapter VI]{Hungerford2003}, which allows a notion of dimension of a field extension, called the \emph{transcendence degree}.

	\begin{Def} \label{def: prereq trans degree}
		Let $F$ be an extension field of $K$. 
		The \concept{transcendence degree} of $F$ over $K$ is the cardinal given by $|S|$, where $S$ is any transcendence base of $F$ over $K$.
	\end{Def}
	
	\begin{Remark}
		Let $F$ be an extension field of $K$ and $S$ a subset of $F$ of cardinality strictly less than the transcendence degree of $F$ over $K$. 
		Then, it follows from \cite[Theorem 1.5 Chapter VI]{Hungerford2003} that there is $u$ in $F$ such that $u$ is a transcendental over $K(S)$.
In other words, if $C$ is an algebraic closure of $K(S)$, 
		then $F\backslash C \neq \emptyset$.

	\end{Remark}

	\subsection{Real analysis}\label{subsection: prereq analysis}

	\begin{Def}\label{def: prereq uniformly continuous}
		Let $f\colon A\subseteq \RR^n \to \RR$. 
		We say that $f$ is \concept{uniformly continuous} on $A$ if 
		\[ \forall \varepsilon>0 \, \exists \delta>0  \, \forall x,y\in A \, \left( \|x-y\|<\delta \implies |f(x)-f(y)|<\varepsilon	\right), \]
		where $\|\cdot\|$ denotes the norm of a vector in $\RR^n$, and $|\cdot|$ denotes the absolute value of a real number.
	\end{Def}

Recall that a continuous function on a compact set is uniformly continuous {\cite[Theorem 4.19]{Rudin1964principles}}, and
  $A\subseteq \RR^n$ is compact iff it is closed and bounded (see for example \cite[Theorem 2.41]{Rudin1964principles}). 
	In particular, a continuous function on the closure of a ball in $\RR^n$ is uniformly continuous. 

\begin{Def}
Let $f\colon A \subseteq \RR^n \to \RR^m$ be a function. We say that $f$ is \concept{Cauchy-continuous} iff 
\[ \forall \{x_l\}_{l<\omega}\subseteq A \text{ Cauchy}, \{f(x_l)\}_{l<\omega} \text{ is Cauchy.} \]
\end{Def}

Any Cauchy-continuous function from $A\subseteq \RR^n$ to $\RR^m$ can be extended to a continuous (and hence Cauchy-continuous) function defined on the Cauchy completion of $A$.
It is a known fact that every uniformly continuous function is also Cauchy-continuous. 
Then, in a compact set, continuous, Cauchy continuous and uniformly continuous are equivalent. This will be used later.

An essential element for the proof of our main theorem is proving a version of the Implicit Function Theorem (Lemma \ref{lemma: implicit function theorem plus} and Lemma \ref{lemma: implicit function th plus new version}). 
So we start by stating the usual one.
	
		\begin{Teo}[Implicit Function Theorem {\cite[Theorem 2-12]{Spivak1965}}] \label{th: prereq implicit function theorem}
Suppose $f\colon \RR^n\times \RR \to \RR$ is continuously differentiable in an open set containing $(a,b)$ and $f(a,b)=0$. 
If $\frac{\partial f}{\partial x_{n+1}} (a,b)\neq 0$, there is an open set $A\subseteq \RR^n$ containing $a$ and an open set $B\subseteq \RR$ containing $b$, with the following property: for each $t\in A$ there is a unique $g(t)\in B$ such that $f(t,g(t))=0$. 
The function $g$ is differentiable. 
		\end{Teo}
		
	\begin{Lema}\label{lemma: IFT rational balls}
		In Theorem \ref{th: prereq implicit function theorem}, one can assume that the open sets $A$ and $B$ are rational balls, namely, each is a ball with center in a rational point and of rational radius. 
	\end{Lema}
	
	\begin{Dem}
		Apply Theorem \ref{th: prereq implicit function theorem} to $f$ and the point $(a,b)$ obtaining open sets $A$ and $B$ such that
		\begin{equation}\label{eq: IFT for A B}
			\forall t\in A \exists! u\in B \text{ such that } f(t,u)=0.
		\end{equation}
		
		Choose rational balls $A'\subseteq A$ and $B'\subseteq B$ such that $a\in A'$ and $b\in B'$.
		Apply the Implicit Function Theorem to $f\restriction A'\times B'$ and the point $(a,b)$, and obtain $A''\subseteq A'$ and $B''\subseteq B'$ such that 
		\begin{equation}\label{eq: IFT for A'' B''}
	\forall t\in A''\, \exists! u\in B'' \text{ such that } f(t,u)=0.
		\end{equation}
		
		Shrink one more time $A''$ to $A'''$ such that $A'''$ is a rational ball and $a\in A'''$. 
		By \ref{eq: IFT for A'' B''}, we get 
		\begin{equation}\label{eq: IFT for A''' B''}
	\forall t\in A''' \,\exists! u\in B'' \text{ such that } f(t,u)=0.
		\end{equation}
		Since $B''\subseteq B'$ and \ref{eq: IFT for A''' B''}, we obtain 
		\begin{equation}\label{eq: IFT for A''' B'}
				\forall t\in A'''\, \exists u\in B' \text{ such that } f(t,u)=0.
		\end{equation}
		From \ref{eq: IFT for A B} and $B'\subseteq B$, we get that for every $t\in A'''$ there is at most one $u\in B'$ such that $f(t,u)=0$. 
		Combining this with \ref{eq: IFT for A''' B'}, we finally get that
		\[ 	\forall t\in A'''\, \exists! u\in B' \text{ such that } f(t,u)=0.\]
		Notice that $A'''$ and $B'$ are rational balls and have the property described in Theorem \ref{th: prereq implicit function theorem}.
		\end{Dem}

	\begin{Def}\label{def: prereq analytic function}
		Let $f\colon U \to  \mathbb{C}$ be a function defined in an open set $U$ of the complex plane $\mathbb{C}$. 
		We say that $f$ is \concept{analytic} in $U$ iff it is representable by power series, namely, for every $a\in U$ there are $r\in \RR_{>0}$ and coefficients $\{c_n\}_{n<\omega}$ such that the following series converges
		\[ \sum_{n=0}^{\infty} c_n (z-a)^n, \]
		and it is equal to $f(z)$
		for all $z$ in the disk $D(a,r)=\{z \colon |z-a|< r\}$.
	\end{Def}
	
Notice that, if $f$ is analytic, then it has derivatives of all orders, and
	\[c_n = \frac{f^{(n)}(a)}{n!}\]
	for every $n<\omega$. 
	
	\begin{Remark}
		Sums, products and compositions of analytic functions are analytic. 
		The exponential function and polynomials are analytic in the whole plane, $\frac{1}{z}$ is analytic in $\mathbb{C}\backslash \{0\}$, each branch of $\log z$ is analytic in $\mathbb{C}\backslash \RR_{\leq 0}$.
	\end{Remark}

	We will also use the following property of analytic functions.
	
	\begin{Teo}[{\cite[Theorem 10.18]{Rudin1987realcomplex}}] \label{th: prereq analytic vanishes in a limit set}
		Let $U$ be an open connected set, and let $f\colon U \to \mathbb{C}$ be an analytic function.
		If $f$ vanishes in a set that has a limit point in $U$, i.e., $f(a)=0$ for every $a$ in that set, then $f\equiv 0$ in $U$.
	\end{Teo}

	\section{The transcendence degree of really closed fields}\label{section: context}
	
	In this section we will approach some algebraic questions related to different subfields of $\RR$, which are actually $\RR$ in some inner model. 
	First, we state a folklore result which is the starting point of this paper. 
	Then we state a generalization in one of the directions mentioned: taking out the condition that the real added is a Cohen real.
	We present the proof given by Ben de Bondt in private conversations. 
	For the main theorem, only Proposition \ref{prop: R really closed implies trans degree c} will be used.

	\begin{Lema}[Folklore] \label{lemma: trascendence degree one cohen}
		Let $M$ be a model of ZFC and $g$ be a $\mathbb{C}$-generic filter over $M$. 
		In $M[g]$, the transcendence degree of $\mathbb{R}$ over the algebraic closure relative to $\mathbb{R}$ of $\mathbb{R} \cap M$ is maximal (i.e. the cardinality of $\mathbb{R}$). 
	\end{Lema}
	
			In particular, if $y$ is a Cohen real over $V[x]$, then in \(V[x,y]\) the transcendence degree of $\mathbb{R}$ over \(\overline{\mathbb{R}\cap V[x]} \) is $\mathfrak{c}$. \\

It is natural to ask for which forcings that add reals Lemma \ref{lemma: trascendence degree one cohen} holds.
	By private (independent) conversations with Ben de Bondt and Elliot Glazer, we know that Lemma \ref{lemma: trascendence degree one cohen} is true for any type of reals added to $V$.
	This result will follow from Proposition \ref{prop: R really closed implies trans degree c}, but before that we need some definitions.
	We will give the proofs of Proposition \ref{prop: R really closed implies trans degree c} and Corollary \ref{coro: trans degree for any inner model with different reals} by Ben De Bondt. 
	
	\begin{Def}
		Let $P \colon \RR_{>0} \to \RR $. 
		We say $P$ is a \concept{generalized polynomial with coefficients} $\{a_i,r_i\}_{i<k}\subseteq \RR$ for some $k<\omega$ iff $r_i\neq r_j$ for $i\neq j$, $a_i\neq 0$ for $i<k$, and 
		\[P(x)\ = \sum_{i<k} a_i x^{r_i}.\]
	\end{Def}
	
	\begin{Def}
		Let $R \subseteq \RR$ be a subfield of the reals. 
		We say $R$ is \concept{really closed} if it is closed under roots of generalized polynomials with coefficients in $R$.  
		Namely, if $P$ is a generalized polynomial with coefficients $\{a_i,r_i\}_{i<k}\subseteq R$ and $x\in \RR_{>0}$ is such that $P(x)=0$,  then $x\in R$.
	\end{Def}

	We will now show that if $P$ is a generalized polynomial then the set of roots $S=\{x \mid P(x)=0\}$ is finite.
	Clearly $S$ is definable over $(\RR, +, \cdot, <, \exp )$ with parameters $\{a_i,r_i\}_{i<k}$. 
	By $o$-minimality \cite{Wilkie1996}, $S$ is a finite union of singletons and intervals. 
	Since $P$ is an analytic function over $\RR_{>0}$, it cannot vanish in an interval, otherwise it would be constantly $0$ (see Theorem \ref{th: prereq analytic vanishes in a limit set}). 
	Therefore $S$ is a finite union of singletons, i.e., a finite set.

	\begin{Prop} \label{prop: R really closed implies trans degree c}
		If $R\subsetneq \RR$ is a really closed subfield of $\RR$, the transcendence degree of $\RR$ over $R$ is continuum. 
	\end{Prop}
	
	\begin{Dem}
		If $|R|< \mathfrak{c}$, the transcendence degree of $ \RR$ over $R$ is $\mathfrak{c}$ by cardinality. 
		Namely, if $B$ is algebraically independent over $R$, and 
		\[S = \{x \mid \exists P \text{ with coefficients in } R \text{ such that }  P(x,{b})=0 \text{ for some } {b}\in B^{<\omega}  \},\] then 
		$|S|=\max\{|B|, |R|\}$. 
		If $B$ is a transcendence base, then $S= \RR$ and therefore we have $|B|= \mathfrak{c}$. 
		
		If $|R|=\mathfrak{c}$, pick $x\in \RR\backslash R$ such that $x>0$, and take $T\subseteq R$ linearly independent over $\mathbb{Q}$ of cardinality $\cc$. 
		Let $B= \{x^t \}_{t\in T}$. 
		We claim that $B$ is algebraically independent over $R$. 
		Suppose not, then there are $k\in \omega$, $\{t_i\}_{i<k}\subseteq T$ and $P$ a non-zero polynomial over $R$ with $k$ variables such that 
		\[P\left(x^{t_0}, \dots, x^{t_{k-1}}\right)=0.\]
		Specifically, there are $J\subseteq \omega^k$ and non-zero coefficients $\{a_j\}_{j\in J} \subseteq R$, such that
		\[\sum_{j\in J} a_j x^{t_0j_0}  \dotsm  x^{t_{k-1}j_{k-1}} =0, \] 
		where for each $j\in J \subseteq \omega^k$ we denote its coordinates by $(j_0, \dots, j_{k-1})$.
		We can rewrite this by grouping the exponents, namely, 
		\[\sum_{j\in J} a_j x^{s_j} =0,\]
		where $s_j = \sum_{i<k} j_i t_i$. Notice that all of these $\{s_j\}_{j \in J}$ are different since $\{t_i\}_{i<k}\subseteq T$ and $T$ is linearly independent over $\mathbb{Q}$. 
		So $x$ is the root of a generalized polynomial with coefficients in $R$. 
		Since $R$ is a really closed subfield of $\RR$, we get that $x\in R$, which is a contradiction. 
		Therefore, $B$ is algebraically independent over $R$ and the transcendence degree of $\RR$ over $R$ is $\mathfrak{c}$. 
	\end{Dem}
	
	\begin{Coro} \label{coro: trans degree for any inner model with different reals}
		If $M$ and $N$ are models of set theory such that $M\subseteq N$ and $\RR^M\subsetneq \RR^N$, then in $N$ the transcendence degree of $\RR^N$ over $\RR^M$ is $(2^{\aleph_0})^N$.
		In particular, if $N$ is a forcing extension of $M$ via a forcing that adds reals, the same conclusion holds. 
	\end{Coro}
	
	\begin{Dem}
		In $N$, $\RR^M$ is a subfield of $\RR^N$. 
		Recall that a generalized polynomial has always finitely many roots. 
		If a generalized polynomial $P$ with coefficients $\{a_i,r_i\}_{i<k}$ in $\RR^M$ has $n$ roots in $N$, this is expressed by the formula
		\[ \exists z_0, \dots, z_{n-1} \in \RR_{>0} \; P(z_0)=\dots=P(z_{n-1})=0 ,\]
		which is $\Sigma^1_1(c)$ for $c=\bigoplus_{i<k}(a_i\oplus r_i)$, and therefore absolute between $M$ and $N$. 
		So the roots of $P$ belong to $M$ and  $\RR^{M}$ is a really closed subfield of $\RR^N$.
		By applying Proposition \ref{prop: R really closed implies trans degree c} in $N$, we get that the transcendence degree of $\RR^N$ over $\RR^M$ is continuum. 
	\end{Dem}
	
	\section{The Main Theorem}\label{section: the theorems}
	
	In this section we present the elements that are necessary for the proof of our main result (Theorem \ref{th: transcendence degree many cohen}). 
	We also include the proof for $|X|=2$ (see Theorem \ref{th: transcendence degree two cohen}), which has all the relevant ideas without the obstacle of notation. 
	We introduce the concept of \emph{$V$-continuous dependence} (Definition \ref{def: depends V continuously on S}), which will be the key element to distinguish the fields $F$ and $\RR^{V[X]}$. 
	
	\begin{Not}
		Let $S\subseteq \RR$. 
		We write $\overline{S}^{\field}$ to denote the minimal subfield of $\RR$ containing $S$. 
		We write $\overline{S}^{\genpol}$ to denote the (real) closure under roots of generalized polynomials with coefficients in $S$. 
		Finally, we write $\overline{S}^{\expfield}$ for the set defined as follows: 
		\[ \overline{S}^{\expfield} := \bigcup_{n<\omega} S_n, \]
		where $S_0:=S$, and $S_{n+1} :=\overline{\overline{S_n}^{\genpol}}^{\field}$.
	\end{Not}\\

	Let $S$ be a subset of $\RR$. 
	Then $\overline{S}^{\expfield}$ is the smallest really closed subfield of $\RR$ containing $S$.
	Clearly, $0,1\in \overline{S}^{\expfield}$. 
	If $a, b, c \in \overline{S}^{\expfield}$ then there is $n\in \omega$ such that $a,b,c \in S_{n+1}$. Since $S_{n+1}$ is a field, $ab^{-1}-c\in S_{n+1}\subseteq \overline{S}^{\expfield}$.  Thus, $\overline{S}^{\expfield}$ is a field. 
	Let $P$ be a generalized polynomial with coefficients on $\overline{S}^{\expfield}$ and let $z\in \RR_{>0}$ be such that $P(z)=0$. 
	Let $\{a_i, r_i\}_{i\in k}\subseteq \overline{S}^{\expfield}$ be the coefficients of $P$. 
	There is $n$ such that $\{a_i, r_i\}_{i\in k}\subseteq S_n$. 
	Then $z\in \overline{S_n}^{\genpol} \subseteq S_{n+1} \subseteq \overline{S}^{\expfield}$, and therefore $\overline{S}^{\expfield}$ is closed under roots of generalized polynomials.

	\begin{Lema}\label{lemma: implicit function theorem plus}
		Let $f \colon \RR^n\times\RR \to \RR$ be a continuously differentiable function in $V$. 
		Let $V[g]$ be a forcing extension of $V$ and $\bar{f}$ the version of $f$ in $V[g]$. 
		Suppose that $a \in \RR^{n}$ and $b \in \RR$ are such that, in $V[g]$, $\bar{f}(a,b)=0$ and $\frac{\partial \bar{f}}{\partial x_{n+1}} |_{(a,b)}\neq 0$.
		
		Then in $V$ there is an open set $U\subseteq \RR^n$ and a continuously differentiable function $h \colon U \to \RR$ such that its version $\bar{h}$ in $V[g]$ satisfies  
		$\bar{h}(a)=b$ and $\bar{f}(t, \bar{h}(t))=0$ for all $t\in \bar{U}$, where $\bar{U}$ is the version of $U$ in $V[g]$, and $a\in \bar{U}$.
	\end{Lema}

	Lemma \ref{lemma: implicit function theorem plus} is just an application of the Implicit Function Theorem. 
	But for our purposes we need to check precisely that in this case the implicit function comes from the ground model $V$.

	\begin{Remark}\label{rem: cont dif is absolute}
		If $f: \RR^{n+1} \to \RR$ is a continuously differentiable function in $V$, then its version $\bar{f}$ in $V[g]$ is also continuously differentiable, as will be proved below.
	\end{Remark}
	
	Let $h: \RR^{n} \to \RR^m$ be a continuous function in $V$. Notice that it is determined by $h\restriction \QQ^n$. 
	We can represent $h$ by a real $r\in \RR^V$ coding $h\restriction \mathbb{Q}^{n}$ (which is a sequence of countably many reals). 
		More precisely, fixing bijections $c$ and $d$ in $V$ so that $c\colon \QQ^{<\omega}\to \omega$ and a $d\colon \RR \to (\RR^{<\omega})^\omega$, we say that ``$r$ codes a continuous function from $\RR^n$ to $\RR^m$''  iff $\varphi(n,m,r)$ holds true, where
	
			\begin{align*}
		\varphi(n,m,r): & \forall (c_0,r_0)\in \QQ^n\times \QQ \; \forall\{q_l\}_{l<\omega} \subseteq \QQ^n \cap \overline{B}(c_0,r_0) \\
&	\{q_l\}_{l<\omega}	\text{ Cauchy } \rightarrow \{d(r)(c(q_l))\}_{l<\omega} \text{ is Cauchy}.
	\end{align*}
	Here $\overline{B}(c_0,r_0)$ denotes the closed ball with center $c_0$ and radius $r_0$.
	Notice that, if $r$ codes a continuous function $h$ then $d(r)(c(q))$ for $q\in \QQ^n$ is just $h(q)$.
	In other words, $r$ is coding what $h\restriction \QQ^n$ will be.
	Also, "being Cauchy" can easily be written as a $\bm{\Delta^1_0}$ formula. 
	Therefore  $\varphi(n,m,r)$ is $\bm{\Pi^1_1}$.
	Finally, $\phi(n,m,r)$ states that $r$ codes a continuous function because continuity and Cauchy continuity are equivalent in compact sets, and continuity is a local property.

Let $f_r$ be a continuous function coded by a real $r$ such that $\varphi(n,m,r)$ holds true. 
Then $f_r$ is unique and it is given by the following formula $\psi(r,x,y)$ that states that $f_r(x)=y$:
\[\psi(r,x,y): \forall \varepsilon \in \QQ_{>0} \exists N \in \omega \forall l \geq N \; \|d(r)(c(\lfloor x \rfloor_l))-y\|<\varepsilon.\]
Here, $\lfloor x \rfloor_l$ denotes the rational number given by truncating $x$ in $l$ digits after the comma.
Let $r$ be a real in $V$ coding a continuous function $f_r$ in $V$. 
Since $\varphi(n,m,r)$ is $\bm{\Pi^1_1}$ formula so absolute between $V$ and $V[g]$, it makes sense to talk about $f_r$ and $\bar{f_r}$ which are given by $\psi(r,\cdot, \cdot)^V$ and  $\psi(r,\cdot, \cdot)^{V[g]}$ respectively. 
Since $\psi(r, x, y)$ is $\bm{\Delta^1_0}$, it is absolute and $\bar{f_r}$ extends $f_r$.

Let $f\colon \RR^n \to \RR$ be a differentiable function and let $g\colon \RR^n \to \RR^n$ be a function. 
We write $g=\nabla f$ iff $g$ is the gradient of $f$, namely, $g(a)=\left( \frac{\partial f}{\partial x_1}(a), \dots , \frac{\partial f}{\partial x_n}(a) \right)$ for every $a\in \RR^n$.
If $r,s\in \RR^V$, let $\phi(n,r,s)$ the following formula stating that $r$ and $s$ code continuous functions $f_r\colon \RR^n \to \RR$ and $f_s\colon \RR^n \to \RR^n$  such that $\nabla f_r =f_s$. 
\begin{align*}
\phi(n,r,s)\colon & \varphi(n,1,r) \wedge \varphi(n,n, s)\wedge \\
							  & \forall x_0 \in \RR^n \forall \varepsilon \in \QQ_{>0} \exists \delta \in \QQ_{>0} \forall x\in \RR^n \\
							  &  \|x-x_0\|<\delta \rightarrow \left| \frac{f_r(x)-f_r(x_0)-\tup{f_s(x_0), x-x_0}}{\|x-x_0\| } \right| <\varepsilon. 
\end{align*}
Notice this is a $\bm{\Pi^1_1}$ formula,  therefore absolute between $V$ and $V[g]$.
All of this is to conclude that the function $\bar{f}$ is also continuously differentiable in $V[g]$.
Finally, if $f$ is continuously differentiable function in $V$ defined in an open set $U$ of $\RR_n$ and its image is contained in $\RR$, then its version $\bar{f}$ in $V[g]$ is also continuously differentiable.
This follows from the fact that differentiability and continuity are local notions and that any open set around a point contains a closed ball around the point and every continuously differentiable function restricted to a closed ball can be extended to a continuously differentiable function defined in all $\RR^n$.

Finally, notice that if $U$ is an open set in $\RR^n$, then it is an union of countably many rational balls. 
If $U$ is an open set in $V$, it can be coded by one real $r$ in $V$, and there is a formula $\Phi(\cdot,r)$ describing the set $U$.
There is a natural interpretation of $U$ in $V[g]$, $\bar{U}$, which is an open set in $V[g]$ and it is given by $\Phi(\cdot, r)^{V[g]}$.
\\

	\begin{Demof}{of Lemma \ref{lemma: implicit function theorem plus}}	By Remark \ref{rem: cont dif is absolute}, $\bar{f}$ is also continuously differentiable in $V[g]$. 
		We can apply then the Implicit Function Theorem in $V[g]$, and we obtain $A$, $B$ open sets in $\RR^n$ and $\RR$ respectively, with $a\in A$ and $b\in B$ such that 
		\begin{equation}\label{eq: lemma IFT}
			\forall t \in A \, \exists ! u\in B \, f(t,u)=0
		\end{equation}
	By Lemma \ref{lemma: IFT rational balls}, we can assume $A$ and $B$ are rational balls, so $A$ and $B$ are the versions in $V[g]$ of $U$ and $W$ rational balls in $V$, respectively. 
	Namely, $A=\bar{U}$ and $B=\bar{W}$.
	The Implicit Function Theorem also states that the function $h \colon A\to B$ given by the formula \ref{eq: lemma IFT} is differentiable.
	Let us re write the formula \ref{eq: lemma IFT} to analyze its level of complexity as follows: 
	\begin{align*}
		\forall t \in \bar{U} \, &\exists u \in \bar{W} \, \psi(r,(t,u),0), \text{ and} \\
		\forall t\in \bar{U} \, &\forall u\in \bar{W} \, \forall u' \in \bar{W} \, \left( \psi(r,(t,u),0)\wedge \psi(r,(t,u'),0) \rightarrow u=u' \right).
	\end{align*}
	
	Here $r$ is the real that codes the function $f$ and $\psi(r,(t,u),0)$ is the formula stating $f(t,u)=0$.
	Notice that since $f$ is in $V$ then $r\in \RR^V$.
	The formula above is $\bm{\Pi^1_2}$. 
	By Shoenfield's absoluteness, we get that the statement \ref{eq: lemma IFT} is also true in $V$: 
	\[V\models 	\forall t \in U \, \exists ! u\in W \, f(t,u)=0 \]
	This defines a function $h_0\in V$. 
	It is clear that $h_0\colon U \to W$ in $V$ is precisely $h\restriction (\RR^n\cap V)$. 
	Since $h$ is continuously differentiable  in $V[g]$, $h_0$ is continuously differentiable in $V$ and $h=\bar{h_0}$.
	\end{Demof}
	
	\begin{Remark}
		From now onwards, we will drop the notation $\bar{f}$ and we will just write $f$ for a continuous function independently of the model in which the function is considered.  
		Similarly, if $U$ is an open set in $\RR^n$ in $V$ we will not write $\bar{U}$ anymore but just $U$ instead, for the version of $U$ in any forcing extension $V[g]$.
	\end{Remark}
	
	\begin{Def}\label{def: depends V continuously on S}
		Working in $V[g]$ (a forcing extension of $V$): 
		Let $S\subseteq \RR$ and $z\in \RR$. 
		We say that $z$ \concept{depends $V$-continuously on} $S$ if there is $F\colon \RR^k \to \RR$ a continuous function in $V$ with $k<\omega$ and $s\in S^{k}$ such that $F(s)=z$.
	\end{Def}
	
	Notice that the interesting case is when $z\not\in S$ and $z\not\in \RR^V$. 
	Otherwise, $F=\text{id}$ or $F\equiv z$ would trivially witness $V$-continuity.
	Observe that $V$-continuous dependence on some set $S$ is a local property. 
	So the function $F$ that witnesses this does not need to be defined in the full space $\RR^k$, but rather in an open set $U\subseteq \RR^k$ such that $s\in U$.\\
	
	Using Definition \ref{def: depends V continuously on S} we can restate Lemma \ref{lemma: implicit function theorem plus} as follows: 
	
	\begin{Lema} \label{lemma: implicit function th plus new version}
		Working in $V[g]$ a forcing extension of $V$:
		Let $f\colon \RR^n \times \RR \to \RR$ a continuously differentiable function in $V$. 
		Let $(a,b)\in \RR^{n}\times \RR$ such that ${f}(a,b)=0$ and $\frac{\partial {f}}{\partial x_{n+1}} |_{(a,b)}\neq 0$.
		Then $b$ depends $V$-continuously on (the set of coordinates of) $a$. 
	\end{Lema}

	Observe that if $z\in \overline{F}$ (the algebraic closure relative to $\RR$ of a subfield $F$),
	by definition there is a polynomial $P$ with coefficients in the field $F$  such that $P(z)=0$. 
	We can assume $P$ has also the property that $P'(z)\neq 0$. 
	This is because a non-zero polynomial cannot vanish in all its derivatives, so we can exchange $P$ for the highest derivative of $P$ which vanishes at $z$, satisfying then the conditions of Lemma \ref{lemma: implicit function th plus new version}.
	As Proposition \ref{prop: gp can be derivative nonzero} shows, this will also be the case for the closure by roots of \emph{generalized} polynomials.
	
	\begin{Prop}\label{prop: gp can be derivative nonzero}
		Let $z\in \overline{F}^{\genpol}$ (the closure of a field $F$ by roots of generalized polynomials). 
	Then there is a generalized polynomial $P$ with coefficients in the field $F$ such that $P(z)=0$ and $P'(z)\neq 0$. 
	\end{Prop}
	
	\begin{Dem}
		By definition of $\overline{F}^{\genpol}$ there is a generalized polynomial $Q$ with coefficients in $F$ such that $Q(z)=0$.
A generalized polynomial is in particular an analytic function on its domain (Definition \ref{def: prereq analytic function}), which in this case includes $\RR_{>0}$. 
If it vanished in all its derivatives at $z$, then the power series of $Q$ at $z$ would have all coefficients equal to zero.
Then, there would be an open set $U$ containing $z$ for which $Q\equiv 0$. 
By Theorem \ref{th: prereq analytic vanishes in a limit set} we get that $Q\equiv 0$ on its domain, this is a contradiction. 
Therefore $Q$ has a non-zero derivative at $z$. 
Let $P$ be the highest derivative of $Q$ that vanishes at $z$. 
 $P$ is itself a generalized polynomial with coefficients in $F$, has $z$ as a root, and its derivative at $z$ is non-zero.
	\end{Dem}

	\begin{Lema}\label{lemma: exp closure of S is V continously dependant over S}
		Let $V[g]$ be a forcing extension of $V$.
		In $V[g]$, let $S\subseteq \RR$ and $z\in \overline{S}^{\expfield}$. 
		Then $z$ depends $V$-continuously on $S$.
	\end{Lema}
	
	\begin{Dem}
		Recall that 
		\[ \overline{S}^{\expfield} := \bigcup_{n<\omega} S_n, \]
		where $S_0:=S$, and $S_{n+1} :=\overline{\overline{S_n}^{\genpol}}^{\field}$.
		
		Let $z\in \overline{S}^{\expfield}$. 
		Let us prove by induction on $n$ that 
		\[z\in S_n \implies z \text{ depends } V\text{-continuously on } S.\]
		
		For $n=0$ and $z\in S_0=S$, it is clear. 
		Suppose the above statement is true for a fixed $n\in \omega$. 
		Let $z\in S_{n+1}=\overline{\overline{S_n}^{\genpol}}^{\field}$. 
		Then there are polynomials $P$ and $Q$ with coefficients in $\mathbb{Q}$ and a finite set $t\subseteq \overline{S_n}^{\genpol}$ such that $\frac{P(t)}{Q(t)}=z$ (See Theorem \ref{th: prereq description field generated by X}). 
		Clearly the function$\frac{P}{Q}$ witnesses that $z$ depends $V$-continuously on $t$. 
		For each $u\in t$ there are $k<\omega$ and a generalized polynomial $P_u$ with coefficients given by some $k$-tuples $a_u$ and $r_u$ in $S_n$, such that $P_{u}(u)=0$. 
		By Proposition \ref{prop: gp can be derivative nonzero}, we can assume that $P_u'(u)\neq 0$. 
		Consider the function $f\colon \RR^{k}\times \RR^k \times \RR_{>0} \to \RR$ defined by 
		\[f(a,r,x) = \sum_{i=0}^{k-1} a_i x^{r_i},\]
for every $a$, $r\in \RR^k$ and $x\in \RR_{>0}$.		
		Then $f$ is continuously differentiable, $f(a_u, r_u, u)=0$ and $\frac{\partial f}{\partial x}(a_u,r_u,u)=P_u'(u)\neq 0$.
		Applying Lemma \ref{lemma: implicit function th plus new version}, $u$ depends $V$-continuously on the set $w_u$ given by the elements of the $k$-tuples $a_u$ and $ r_u$. 
		
		We know that $z$ depends $V$-continuously on $t$ and each element $u$ of $t$ depends $V$-continuously on $w_u$. 
		Then $z$ depends $V$-continuously on the finite set $w:=\bigcup_{u\in t} w_u \subseteq S_n$, witnessed by the composition of the functions that witness each step. 
		By inductive hypothesis, each element of $w$ depends $V$-continuously on $S$.
		Therefore $z$ depends $V$-continuously on $S$. 
	\end{Dem}
	
	\begin{Def}
		Let $f \colon C\subseteq \RR^k \to \RR$ be a continuous function defined on a compact set $C$.
		We say that $\delta \colon \mathbb{Q}_{>0} \to \mathbb{Q}_{>0}$ \concept{witnesses the uniform continuity of} $f$  if 
		\[ \forall\varepsilon\in \mathbb{Q}_{>0} \forall x, x'\in C \; \|x-x'\|<\delta(\varepsilon) \rightarrow  \left| f(x) - f(x')\right| <\varepsilon. \]
	\end{Def}
	
\begin{Remark}\label{rem: witness for unif cont is absolute}
	Notice this is a $\bm{\Pi^1_1}$ formula therefore absolute between transitive models containing its parameters. 
Recall that $f$ and $\delta$ can be coded by a real. 
\end{Remark}
	
	Now we are ready to prove the main theorem. 
	For better readability, we first include the proof of it when $|X|=2$ (see Lemma \ref{lemma: there is z not V-continuously dependant 2 cohen} and Theorem \ref{th: transcendence degree two cohen}), which simplifies the notation but shows all the relevant ideas of the general proof. 
	
	\begin{Lema} \label{lemma: there is z not V-continuously dependant 2 cohen}
		Let $x$ and $y$ be mutually Cohen-generic reals over $V$. 
		In $V[x,y]$, consider the set $S= \RR^{V[x]} \cup \RR^{V[y]}$.
		Then there is $z\in \RR$ that does not depend $V$-continuously on $S$.
	\end{Lema}
	
	\begin{Dem}
		We will define $z$ from $x$ and $y$.
		For this purpose we need a function $c\colon {^{\omega}}\omega \to [0,1]$ that will help us translate properties of the Baire space to the real line in a nice way\footnote{At a first glance, it seems that is not so relevant which is the $c$ we choose. 	But, for example, a representation of reals as binary sequences does not work for the last part of the proof.}.
		
		Let $a=(a_n)_{n<\omega}$ be an element of ${^\omega}\omega$. 
		Then we define $c(a)\in [0,1]$ as 
		\[c(a) = \sum_{n<\omega} \frac{a_n \bmod 3}{3^{n+1}},  \] 
		where $m \bmod 3 $ denotes the unique number in $\{0,1,2\}$ that has the same residue as $m\in \omega$ in the division by 3. 
		
		Let $s$ be a sequence of natural numbers of length $n$. 
		Then we know that there is a closed interval $I_s$ of length $\frac{1}{3^n}$ such that for any $a\in {^\omega}\omega$: $c(a)\in \text{int}(I_s)$ implies $s$ is an initial segment of $a$; and if $s$ is an initial segment of $a$, then $c(a)\in I_s$.
		
		Think of $x$ and $y$ as elements of the Baire space. 
		Consider $y\circ x\in {^{\omega}}\omega$ (composition of functions) and $z=c(y\circ x)$. 
		By definition of $c$, $z\in \RR^{V[x,y]}$.
		\\
		
		Suppose on the contrary that $z$ depends $V$-continuously on $S$. 
		Then there are $k, l\in \omega$, $U$ an open set in $\RR^{k+l}$, $F\colon U \to\RR$ a continuous function in $V$, and sets $r\subseteq \RR^{V[x]} $, $s\subseteq \RR^{V[y]}$ of sizes $k$ and $l$ such that $F(r,s)=z$.
		Restrict $F$ to a closed ball $C$ in $V$ such that $(r,s)$ belongs to (the version in $V[x,y]$ of) $C$. 
		Recall that $C$ is compact, and therefore $F$ is uniformly continuous in $C$.
		Let $\delta \colon \mathbb{Q}_{>0} \to \mathbb{Q}_{>0}$ be a witness for the uniform continuity of $F$ restricted to $C$ in $V$. 
		Notice that $\delta$ can be coded by a real in $V$, and by Remark \ref{rem: witness for unif cont is absolute} it is also a witness for the uniform continuity of (the version in $V[x,y]$ of) $F$ restricted to (the version in $V[x,y]$ of) $C$ in $V[x,y]$.
		\\
		
		Given a name $\sigma \in V^\mathbf{C}$, we define $\sigma^{\ast}$ and $^\ast \sigma$ in $V^\mathbf{C\times C}$ recursively:
		\begin{align*}
			\sigma^\ast &= \{ (\pi^\ast,(p,q)) \mid (p,q)\in \mathbf{C\times C} \text{ and } (\pi, p)\in \sigma\}, \text{ and } \\
			^\ast\sigma &= \{ (^\ast\pi,(p,q)) \mid (p,q)\in \mathbf{C\times C} \text{ and } (\pi, q)\in \sigma \}.
		\end{align*}
		
		Let $\tau\in V^{\mathbf{C\times C}}$ a name for $z$, $\pi\in V^{\mathbf{C}}$ a name for $r$ and $\sigma\in V^{\mathbf{C}}$ a name for $s$. 
		Let  $\dot{r}=\pi^\ast$ and $\dot{s}={^\ast}\sigma$.
		Then $\dot{r}$ and $\dot{s}$ are $\mathbf{C\times C}$-names for $r$ and $s$, and they only depend on the first and second coordinate respectively.
		
		Let $(p,q)$ be a condition in $\mathbf{C}\times \mathbf{C}$ such that 
		
		\[(p,q) \forc{\mathbf{C}\times \mathbf{C}}{V} 
		(\dot{r},\dot{s})\in \check{C}
		\text{ and } \check{F}(\dot{r},\dot{s})=\tau. \]
		
		Here and in the rest of the proof, we are abusing the check notation to emphasize which objects are defined in the ground model $V$. 
		In particular, any closed ball or open set $C$ and any continuous function $F$ in $V$ can be coded by reals $c$ and $f\in \RR^V$, as we discussed in Remark \ref{rem: cont dif is absolute}.
	When we write ``$r\in C$'', ``$\dot{r}\in \check{C}$'',``$F(r)=z$'', and ``$\check{F}(\dot{r})=\tau$'', we actually mean $\Phi(r,c)$, $\Phi(\dot{r},\check{c})$, $\psi(f,r,z)$, and $\psi(\check{f},\dot{r},\tau)$ respectively, where $\Phi$ and $\psi$ are the formulas describing $C$ and $F$ from the reals $c$ and $f$ coding them.

		We will find a condition $(\tilde{p},\tilde{q})$ below $(p,q)$ that forces the opposite statement, reaching a contradiction.
		
		First let us assume $\lh(p)=\lh(q)=m$. 
		We will extend $(p,q)$ in several steps until getting to the desired $(\tilde{p}, \tilde{q})$, as Figure \ref{fig:transdegree2cohen-step34} shows.
		
		\begin{enumerate}
			
			\item  Extend $(p,q)$ to $(p, q')$ so that $\lh(q')=\max( \im p)+1$.
			Recall that $x$ and $y$ are the generic associated to the first and second coordinate respectively, thus $(p,q')$ decides $(y\circ x)\restriction m$ but does not decide $(y\circ x)(m)$.
			
			This implies that $(p,q')$ forces that $\tau \in \check{I}$, where $I$ is some interval given by the code $c$. 
			More precisely, $I=I_t$ where $t$ is the sequence given by 
			$(q'_{p_0}, \dots, q'_{p_{m-1}})$.
			
			Let $\varepsilon= \frac{\lh(I)}{6}$ and $\delta=\delta(\varepsilon)$ (this is computed in $V$). 
			
			\item Extend $(p,q')$ to $(p, q'')$ so that it decides $\dot{s} \in \check{J}_\delta$ (some  ball in $\RR^l$  with a rational center and of radius less than $\delta / 4$). 
			
			\item Extend $(p,q'')$ to $(p', q'')$ so that $p'(m)=\dom(q'')$ and $\dom(p')=m+1$ (see Figure \ref{fig:transdegree2cohen-step34}).
			This implies $q''(p'(m))$ is not defined and therefore $(p',q'')$ does not decide the value of $(y\circ x)(m)$.

			\begin{figure}
				\centering
				\includegraphics[width=0.7\linewidth]{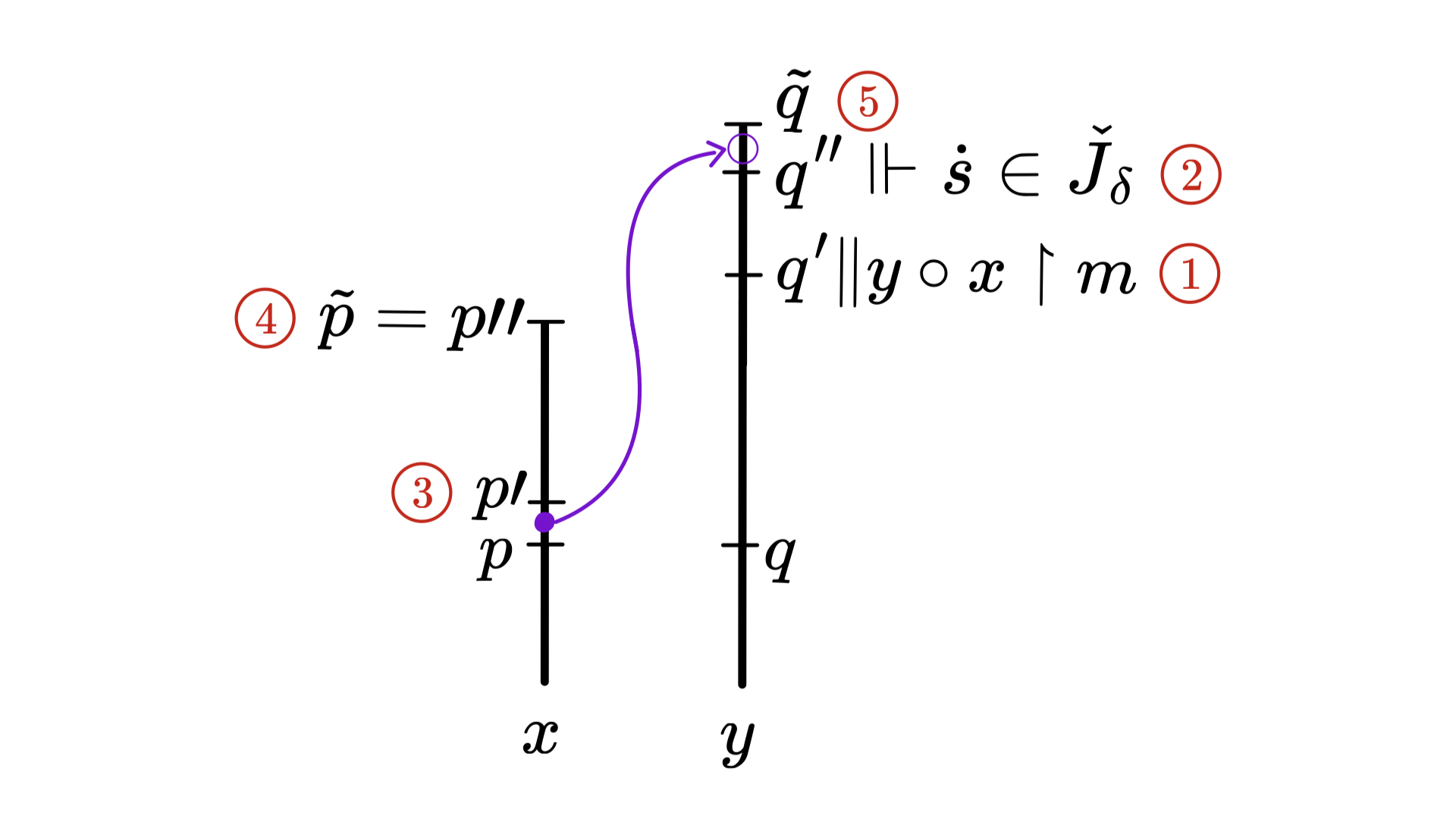}
				\caption{Steps 1--5. Here the arrow represents the fact that $p'(m)=\dom(q'')$, where $m=\lh(p)$.}
				\label{fig:transdegree2cohen-step34}
			\end{figure}
			
			\item Extend $(p',q'')$ to $(p'', q'')$ so that $p'$ decides $\dot{r}\in \check{I}_{\delta}$ (another ball in $\RR^k$ of rational center and radius less than $\delta / 4$).
			
			In $V$, take $(u,v)\in I_{\delta}\times J_{\delta}$.  
			Notice that 
			\[(p'', q'')\forc{\mathbb{C}\times \mathbb{C}}{V}  (\dot{r},\dot{s}), (\check{u}, \check{v}) \in \check{I}_{\delta} \times \check{J}_{\delta} \subseteq \check{C}.\]
			
			If $(u,v)$ and $(r,s)$ are elements of $I_\delta \times J_\delta$, then
			$\|r-u\|<2 \delta/4 = \delta/2$ and $\|s -v\|<\delta/2$.
			Then, $\|(r,s) - (u,v)\| <\delta$. So we obtain 
			\[(p'', q'')\forc{\mathbb{C}\times \mathbb{C}}{V}  \| (\dot{r},\dot{s}) - (\check{u}, \check{v}) \| < \delta.\]
			
			Let $w=F(u,v)\in \RR^V$. 
			By the absoluteness mentioned in Remark \ref{rem: witness for unif cont is absolute},
			\[\text{$\mathbbm{1}$} \forc{\mathbb{C}\times \mathbb{C}}{V} \check{\delta} \text{ witnesses uniform continuity of } \check{F} \text{ on } \check{C}.\]
			
			Then we get that 
			\begin{align}
				(p'', q'')\forc{\mathbb{C}\times \mathbb{C}}{V} &\, |\check{F}(\dot{r}, \dot{s})-\check{F}(\check{u},\check{v})| <\varepsilon, \text{ i.e.,} \nonumber \\
				(p'', q'')\forc{\mathbb{C}\times \mathbb{C}}{V} &\, |\tau-\check{w}|< \varepsilon. \label{eq: transdegree 2 cohen}
			\end{align}
			
			\item Finally, extend $(p'', q'')$ to $(\tilde{p}, \tilde{q})$ so that $\tilde{p}=p''$ and $\lh(\tilde{q})=\lh(q'')+1$ as follows.
			Let $I'=(w-\varepsilon, w+\varepsilon)$. 
			We know that $\lh(I')=2\varepsilon=\frac{\lh I}{3}$, 
			so $I'\cap I \subsetneq I$.
			Notice that $I$ and $I'$ are intervals in $V$, in the sense of their endpoints being reals in $V$.
			Recall that $I$ is an interval $I_t$ for some finite sequence $t$ and $I_t=I_{t{}^\frown0}{\cup} I_{t{}^\frown1} {\cup} I_{t{}^\frown2}$ where each of the intervals $I_{t{}^\frown i}$ is of length $\frac{\lh(I)}{3}$. 
			Since $\lh(I')=\frac{\lh I}{3}$ as well, there is a number $j\in \{0,1,2\}$ such that $I_{t{}^\frown j}\cap I'=\emptyset$.
			Let us assign to $\tilde{q}(\tilde{p}(m))$ the number $j$, so that $z$ avoids being in $I'$.
			
			Then we get 
			\[(\tilde{p}, \tilde{q})\forc{\mathbb{C}\times \mathbb{C}}{V} \tau \not \in \check{I'}, \]
			which contradicts Equation \ref{eq: transdegree 2 cohen}.
		\end{enumerate}
		Thus, $z$ as we defined it does not depend $V$-continuously on $\RR^{V[x]}\cup \RR^{V[y]}$.
	\end{Dem}
	
	\begin{Remark}
		The composition of two mutually generic Cohen reals is a Cohen real. 
	\end{Remark}
	
	\begin{Teo}\label{th: transcendence degree two cohen}
		Let $x$ and $y$ be mutually generic Cohen reals over $V$. 
		In $V[x,y]$, consider the minimum field $F\subseteq \RR$ such that $F\supseteq \RR^{V[x]} \cup \RR^{V[y]}$.
		Then, in $V[x,y]$ the transcendence degree of $\RR$ over ${F}$ is continuum. 
	\end{Teo}
	
	\begin{Dem}
		Work in $V[x,y]$.
		Let $S=\RR^{V[x]} \cup \RR^{V[y]}$. 
		Notice that $\overline{S}^{\expfield}\supseteq \overline{F}$, where $\overline{F}$ denotes the real-algebraic closure of $F$.
		By Lemma \ref{lemma: there is z not V-continuously dependant 2 cohen}, there is a real $z$ that does not depend $V$-continuously on $S$. Applying Lemma \ref{lemma: exp closure of S is V continously dependant over S}, we deduce that $z\not\in \overline{S}^{\expfield}$. 
		Recall that $\overline{S}^{\expfield}$ is a really closed subfield of $\RR^{V[x,y]}$. 
		Using Proposition \ref{prop: R really closed implies trans degree c}, we get that the transcendence degree of $\RR^{V[x,y]}$ over $\overline{S}^{\expfield}$ is continuum. 
		Therefore the transcendence degree of  $\RR^{V[x,y]}$ over $\overline{F}$ is also continuum. 
	\end{Dem}
	
	Now we will prove the version of Theorem \ref{th: transcendence degree two cohen} in which we allow to have more than two Cohen reals. 
	For this purpose, we also need a more general version of Lemma \ref{lemma: there is z not V-continuously dependant 2 cohen}, and that is the role of Lemma \ref{lemma: there is z not V-continuously dependant many cohen}.
	
	\begin{Lema} \label{lemma: there is z not V-continuously dependant many cohen}
		Let $X$ be a finite set of mutually generic Cohen reals over $V$. 
		In $V[X]$, consider the set $S= \bigcup_{Y\subsetneq X} \RR^{V[Y]}$.
		Then there is $z\in \RR$ that does not depend $V$-continuously on $S$.
	\end{Lema}
	
	\begin{Dem}
		We will define $z$ from $X$ in a similar way to the one in Lemma \ref{lemma: there is z not V-continuously dependant 2 cohen}.
		Let $c\colon {^{\omega}}\omega \to [0,1]$ the same function defined in the proof of Lemma \ref{lemma: there is z not V-continuously dependant 2 cohen}.
		
		Think of $X$ as a finite subset of the Baire space, $X=\{x_0, \dots, x_{k-1}\}$ and $k<\omega$.
		Consider $x_{k-1}\circ \dotsb \circ x_0 \in {^{\omega}}\omega$ (composition of functions) and let $z=c(x_{k-1}\circ \dotsb \circ x_0)$. 
		
		Suppose on the contrary that $z$ depends $V$-continuously on $S$. 
		Then there are $n \in \omega$, $U\subseteq\RR^{n}$ open, $F\colon U \to\RR$ continuous function in $V$ and $r\subseteq S $ such that $F(r)=z$.
		Restrict $F$ to a closed ball $C$ in $V$ such that $r$ belongs to (the version in $V[X]$ of) $C$. 
		Since $C$ is compact, $F$ is uniformly continuous in $C$.
		Let $\delta \colon \mathbb{Q}_{>0} \to \mathbb{Q}_{>0}$ be a witness for the uniform continuity of $F$ restricted to $C$.
		Notice that $\delta \in V$ and by Remark \ref{rem: witness for unif cont is absolute} it is also a witness for the uniform continuity of (the version in $V[X]$ of) $F$ (the version in $V[X]$ of) restricted to $C$ in $V[X]$.
		
		To simplify notation we will write $X^\ast_i$ for $X\backslash\{x_i\}$. Notice that $S=\bigcup_{i< k} \RR^{V[X^\ast_i]}$.
	Assume without loss of generality that the coordinates of $r\in \RR^n$ are ordered with respect to this covering of $S$, namely, there are $\{r_i\}_{i<k}$ such that $r_i\in \RR^{<\omega}\cap V[X^\ast_i]$ for all $i<k$,  $\sum_{i<k} \lh (r_i) =n$, and $r=(r_0, \dots, r_{k})$.
		
		Given a name $\sigma \in V^\mathbf{C^{k-1}}$, we define $\sigma^{\ast}_i$ in $V^\mathbf{C^k}$ for $i<k$ as follows:
		\[\sigma^\ast_i = \{(\pi^\ast_i,p) \mid p\in \mathbf{C^k}  \text{ and } (\pi, (p_0, \dots, p_{i-1}, p_{i+1}, \dots, p_{k-1}))\in \sigma\}.\]
		
		Let $\tau\in  V^\mathbf{C^{k}}$ be a name for $z$, and for each $i<k$, let $\dot{r_i}\in  V^\mathbf{C^{k}}$ be an adequate name for $r_i$ that only depends on the coordinates $k\backslash\{i\}$, namely, $\dot{r_i}$ is of the form $\sigma^\ast_i$ for some $\sigma\in  V^\mathbf{C^{k-1}}$ name for $r_i$. 
		Let us write $\dot{r}$ for the name $(\dot{r}_0,\dots, \dot{r}_{k-1})$.
		
		Let $q$ be a condition in $ \mathbf{C}^k$ such that 
		\[ q \forc{\mathbf{C}^{k}}{V} 
		\dot{r} \in \check{C}
		\text{ and } \check{F}(\dot{r})=\tau. \]
		
		We will find a condition $\tilde{q}$ below $q$ that forces a contradictory statement, reaching a contradiction.
		Here and in the rest of the proof, similarly to the proof of Lemma \ref{lemma: there is z not V-continuously dependant 2 cohen}, we are abusing the check notation to emphasize which objects are defined in the ground model $V$. 
		In particular, any closed ball or open set $C$ and any continuous function $F$ in $V$ can be coded by reals $c$ and $f\in \RR^V$, as we discussed in Remark \ref{rem: cont dif is absolute}.
		When we write ``$r\in C$'', ``$\dot{r}\in \check{C}$'',``$F(r)=z$'', and ``$\check{F}(\dot{r})=\tau$'', we actually mean $\Phi(r,c)$, $\Phi(\dot{r},\check{c})$, $\psi(f,r,z)$, and $\psi(\check{f},\dot{r},\tau)$ respectively, where $\Phi$ and $\psi$ are the formulas describing $C$ and $F$ from the reals $c$ and $f$ coding them. 
		
		First let us assume $\lh(q_i)=m$ for all $i<k$.  
		Extend $q$ to $p=(p_0, \dots, p_{k-1})$ so that $p_{k-1}\circ \dots \circ p_0\restriction m $ is defined, setting $p_0=q_0$ and recursively  $\lh(p_{i+1})=\max( \im p_{i})+1$ for $i< k-1$.
		Thus, it is forced by $p$ that $\tau\in \check{I}$, where $I$ is some interval given by the code $c$. 
		More precisely, $I=I_t$ where $t$ is $(p_{k-1}\circ \dots \circ p_0 (0), \dots, p_{k-1}\circ \dots \circ p_0(m-1))$.
		Let $\varepsilon= \frac{\lh(I)}{6}$ and let $\delta=\delta(\varepsilon)$, which is computed in $V$. 
		
		Make $ p\backslash\{p_0\}$ also decide $\dot{r}_0\in \check{I}_0$, where $I_0$ is a ball of radius less than $\frac{\delta}{2k}$.
		Notice that $p_{k-1}\circ \dots \circ p_0 (m)$ is not defined yet because we did not extend $p_0$.
		We will extend sequentially $p$ to conditions $p^0, \dots , p^{k-1}$ in $\mathbf{C}^{k-1}$ as Figure \ref{fig:transdegreemanycohen-steps} shows.

		\begin{figure}
			\centering
			\includegraphics[width=1.1\linewidth]{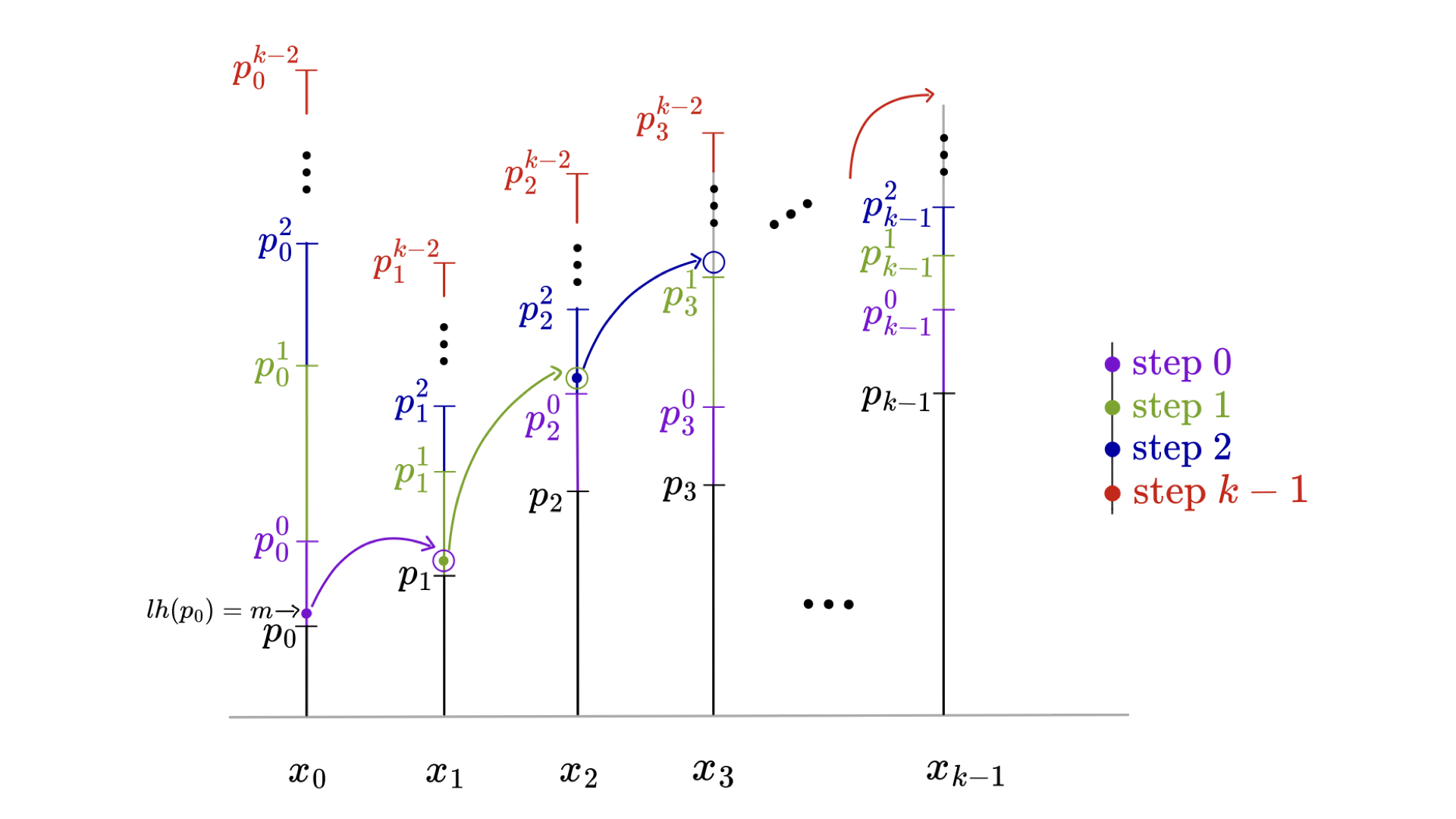}
			\caption{Steps $0$ to $(k-1)$ of the construction. The first arrow represents that in Step 0 we take $p_0^0(m)=\dom(p_1)$, where $m=\lh(p_0)$, and similarly for the other arrows.}
			\label{fig:transdegreemanycohen-steps}
		\end{figure}
		
		We denote coordinates with subindices and the different steps with superindices. 
		\begin{description}
			\item[Step 0.] Define $p^0\in \mathbf{C}^{k}$ extending $p$ such that 
			\begin{itemize}
				\item $p^0_0(m)=\dom(p_1)$,
				\item $p^0_1=p_1$, and
				\item $p^0 \backslash \{p^0_1\}$ decides $\dot{r}_1\in \check{I}_1$, 
				where $I_1$ is a ball in $V$ of radius less than $\frac{\delta}{2k}$.
			\end{itemize}
			\item[Step 1.] Define $p^1\in \mathbf{C}^{k}$ extending $p_0$ such that 
			\begin{itemize}
				\item $p^1_1\circ p^0_0 (m)=\dom(p^0_2)$,
				\item $p^1_2=p^0_2$, and
				\item $p^1 \backslash \{p^1_2\}$ decides $\dot{r}_2\in \check{I}_2$, 
				where $I_2$ is a ball in $V$ of radius less than $\frac{\delta}{2k}$.
			\end{itemize}
		\end{description}
		In general, for $i=1,\dots, k-2$: 
		\begin{description}
			\item[Step i.] Extend $p^{i-1}\in \mathbf{C}^{k}$ to $p^i$ so that 
			\begin{itemize}
				\item $p^i_i\circ \dots \circ p^0_0 (m) = \dom(p^{i-1}_{i+1})$, 
				\item $p^i_{i+1}=p^{i-1}_{i+1}$, and
				\item $p^i \backslash \{p^i_{i+1}\}$ decides $\dot{r}_{i+1}\in \check{I}_{i+1}$, 
				where $I_{i+1}$ is a ball in $V$ of radius less than $\frac{\delta}{2k}$.
			\end{itemize}
		\end{description}
		
		In $V$, take $s\in \prod_{i<k} I_i$, and let $y=F(s)$.
		Notice that:
		
		\begin{align}
			p^{k-2} & \forc{\mathbf{C}^k}{V} \dot{r}, \check{s} \in \prod_{i<k} \check{I}_i, \nonumber\\
			p^{k-2} & \forc{\mathbf{C}^k}{V} \| \dot{r} - \check{s} \|  <\delta \nonumber, \\
			p^{k-2} & \forc{\mathbf{C}^k}{V} \left| \check{F}(\dot{r}) - \check{F}(\check{s}) \right| < \varepsilon, \nonumber \\
			p^{k-2} & \forc{\mathbf{C}^k}{V} \tau \in \check{I'}. \label{eq: transdegree many cohen}
		\end{align}
		
		Here, $I':=(y-\varepsilon, y+\varepsilon)$. 
		Notice that $I'$ is an interval with end points in $V$ and $\lh(I')=2\varepsilon= \frac{\lh(I)}{3}$.
		
		\begin{description}
			\item[Step k-1.] Extend $p^{k-2}\in \mathbf{C}^{k}$ to $p^{k-1}$ such that 
			$p^{k-1}_{k-1}\circ \dots \circ p^0_0 (m)$ is a number $l$ that makes $z$ avoid $I'$, 
			namely, there is $l\in \{0,1,2\}$ such that 
			$I_{t{}^\frown l} \cap I'=\emptyset$.
			This is possible because $\lh I'=\lh I_{t{}^\frown 0}=\lh(I_{t{}^\frown 1})=\lh(I_{t{}^\frown 2}) = \frac{\lh(I)}{3}$.
		\end{description}
		
		Define $\tilde{q}=p^{k-1}$. 
		Then we get 
		\[	\tilde{q}  \forc{\mathbf{C}^k}{V} \tau \not \in \check{I'}, \]
		which contradicts Equation \ref{eq: transdegree many cohen}.
		
		Therefore, $z$ does not depend $V$-continuously on $S=\bigcup_{Y\subsetneq X} \RR^{V[Y]}$.
	\end{Dem}

	\begin{Teo}\label{th: transcendence degree many cohen}
		Let $X$ be a finite set of mutually generic Cohen reals over $V$.  
		In $V[X]$, consider the minimum field $F\subseteq \RR$ such that $F\supseteq \bigcup_{Y\subsetneq X} \RR^{V[Y]}$. 
		Then, in $V[X]$ the transcendence degree of $\RR$ with respect to ${F}$ is continuum. 
	\end{Teo}
	
	\begin{Dem}
		Work in $V[X]$.
		Let $S=\bigcup_{Y\subsetneq X} \RR^{V[Y]}$. 
		Notice that $\overline{S}^{\expfield}\supseteq \overline{F}$, 
		where $\overline{F}$ denotes the real-algebraic closure of $F$. 
		By Lemma \ref{lemma: there is z not V-continuously dependant many cohen}, there is a real $z$ that does not depend $V$-continuously on $S$. 
		Applying Lemma \ref{lemma: exp closure of S is V continously dependant over S}, we deduce that $z\not\in \overline{S}^{\expfield}$. 
		Recall that $\overline{S}^{\expfield}$ is a really closed subfield of $\RR^{V[X]}$. 
		Using Proposition \ref{prop: R really closed implies trans degree c}, we get that the transcendence degree of $\RR^{V[X]}$ over $\overline{S}^{\expfield}$ is continuum. 
		Therefore the transcendence degree of  $\RR^{V[X]}$ over $\overline{F}$ is also continuum. 
	\end{Dem} 
	
	\begin{question}
		Is Theorem \ref{th: transcendence degree many cohen} true for other forcings that add reals? 
	\end{question}

	Theorem \ref{th: transcendence degree two cohen} shows that a transcendence base $B$ of $\RR^{V[x,y]}$ over $F$ (the minimum subfield of $\RR^{V[x,y]}$ such that $F\supseteq \RR^{V[x]} \cup \RR^{V[y]}$) has size continuum. 
	If one would like to extend $B$  to a transcendence base of $\RR^{V[x,y]}$ over $\RR^V$ instead, we need to produce a transcendence base $C$ of $F$ over  $\RR^V$. 
	It turns out that there is a natural candidate for $C$: a transcendence base of $\RR^{V[x]}$ over $\RR^V$ union a transcendence base of $\RR^{V[y]}$ over $\RR^V$, as the following proposition shows.

	\begin{Prop}\label{prop: alg fact}
		Let $x, y$ be Cohen-mutually generic filters over $V$, where $V$ is a model of $\mathsf{ZFC}$.
		Let \(B \subseteq \mathbb{R}^{V[x]}\) be an algebraically independent set over $\mathbb{R}^V$, then $B$ is also algebraically independent over \(\mathbb{R}^{V[y]}\). 
	\end{Prop}
	
	\begin{figure}[ht]
		\centering
		\includegraphics[width=0.8\linewidth]{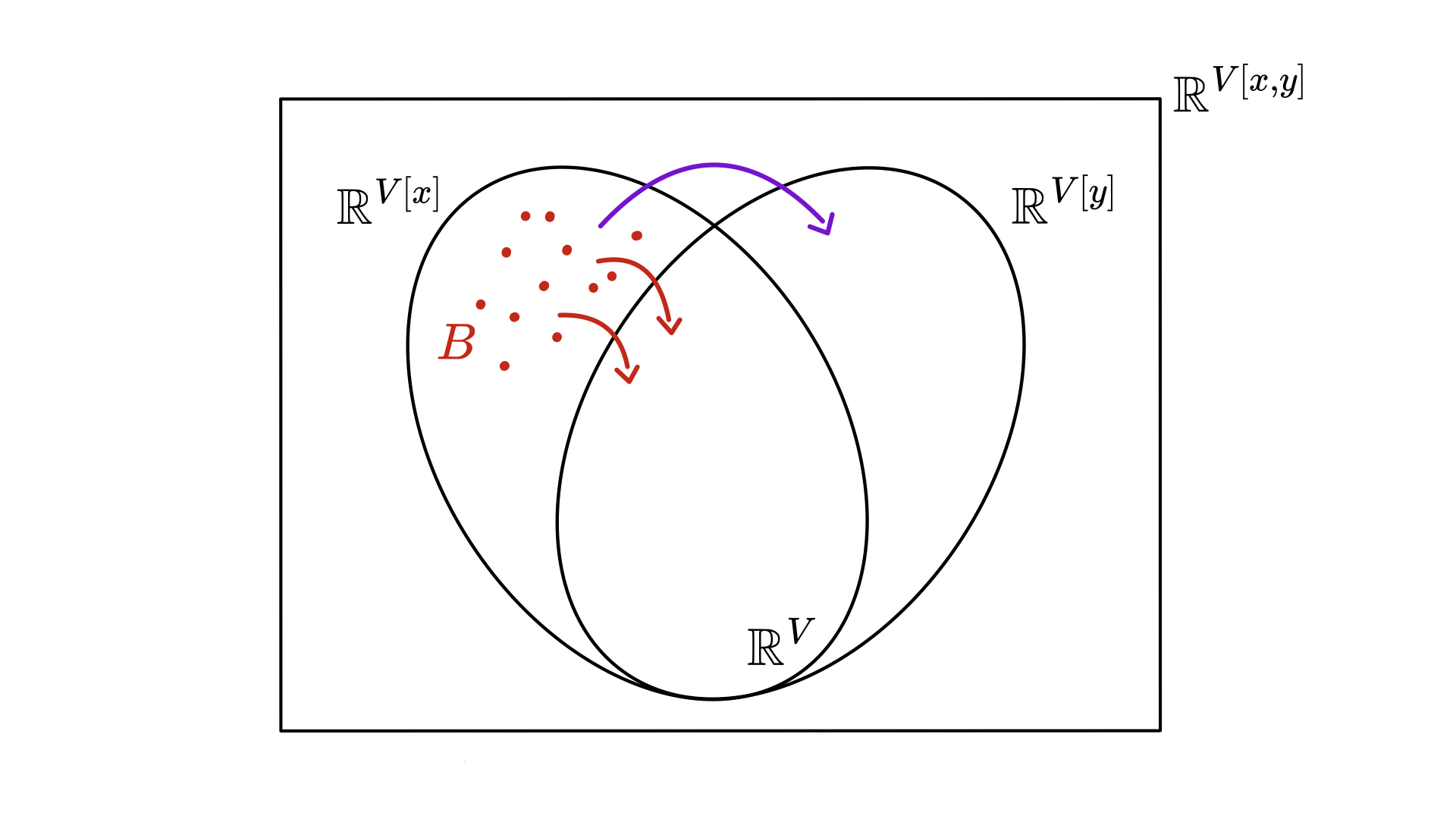}
		\caption{If \(B \subseteq \mathbb{R}^{V[x]}\) is algebraically independent over $\mathbb{R}^V$, then $B$ is also algebraically independent over \(\mathbb{R}^{V[y]}\). }
		\label{fig:algebraiclemma}
	\end{figure}

	\begin{Dem}
		Suppose not, i.e. $B$ is algebraically dependent over \(\mathbb{R}^{V[y]} \) and without loss of generality assume $B$ is finite. 
		Let $n$ be the cardinality of $B$. 
		This means there are some finite multi index set \( \mathcal{J} \subseteq \prescript{n}{}{\omega} \) and some non-zero coefficients \( \{d_J\}_{J\in \mathcal{J}} \subseteq \mathbb{R}^{V[y]} \) such that
		
		\begin{equation} \label{eq: alg fact}
			\sum_{J\in \mathcal{J}} d_J B^J=0,
		\end{equation}
		
		where \(J=(j_0,\ldots, j_{n-1}) \), \(B= \{b_0, \ldots, b_{n-1}\}\) and \( B^J \vcentcolon b_0^{j_0}\cdots b_{n-1}^{j_{n-1}}  \). 
		\\
		
		On the other hand, we can think of \(\mathbb{R}^{V[y]}[B] \) (the minimal ring containing $\mathbb{R}^{V[y]}$ and $B$) as a vector space over the field $\mathbb{R}^{V[y]}$.
		In this context, it is clear that \( \{B^J\}_{J\in{}^n\omega } \) spans $\mathbb{R}^{V[y]}(B)$.
		Therefore, there is \(\mathcal{S} \subseteq \{B^J\}_{J\in {^n}\omega} \) such that $\mathcal{S}$ is a basis of $\mathbb{R}^{V[y]}[B]$ as $\mathbb{R}^{V[y]}$-vector space. 
		Let $\mathcal{I}$ be the corresponding index set, i.e, \(\mathcal{I}\subseteq {^n}\omega \) and \( \mathcal{S}= \{B^J\}_{J \in \mathcal{I}} \). 
		
		Note that \(\mathcal{J}\not \subseteq \mathcal{I}\). 
		$\mathcal{J}$ cannot be a subset of $\mathcal{I}$, since the equation \ref{eq: alg fact} shows a linear dependence of \( \{B^J\}_{J\in \mathcal{J}}\) over $\mathbb{R}^{V[y]}$, and \( \mathcal{S}= \{B^J\}_{J \in \mathcal{I}} \) is a basis, in particular, it is linearly independent over $\mathbb{R}^{V[y]}$. 
		
		Then, there must be an \(I\in \mathcal{J}\backslash \mathcal{I}\), or equivalently, there is an $I \in \mathcal{J}$ such that \(B^I \not \in \mathcal{S}\).
		Now, since \(B^I \in \mathbb{R}^{V[y]}[B]  \), 
		there are unique coefficients \( \{c_i\}_{i=0}^{m-1}\) in $\mathbb{R}^{V[y]}$ and vectors \( \{S_i\}_{i=0}^{m-1}\) in $\mathcal{S}$ such that 
		\begin{equation}
			\sum_{i=0}^{m-1} c_i S_i = B^I. 
		\end{equation}
		
		Since we know 
		\begin{equation}
			V[x, y] \models ``\sum_{i=0}^{m-1} c_i S_i = B^I, \text{ with } \{c_i\}_{i=0}^{m-1} \subseteq \mathbb{R}^{V[y]} \text{''}
		\end{equation}
		
		there is some condition $p\in y \subseteq \mathbf{C}$ such that 
		\begin{equation}
			p \forc{\mathbf{C}}{V[x]} 
			``\sum_{i=0}^{m-1} \tau_i  \check{S_i} = \check{B^I}, \text{ with } \{\tau_i\}_{i=0}^{m-1}\subseteq \mathbb{R}^{V[\dot{g}]}, \text{''}
		\end{equation}
		
		where $\dot{g}$ is the usual name for the generic filter, and $\tau_i$ is a name for $c_i$, i.e., \( \tau_i^{y}=c_i\) for $i=0, \dots , m-1$. 
		Note also that $S_i \in \mathcal{S}= \{B^J\}_{J\in \mathcal{I}} \subseteq V[x]$, which justifies the ``check'' on $S_i$, for $i=0,\dots , m-1$. 
		
		Now, let split $y$ into two mutually Cohen generics $y_1, y_2$ over $V[x]$ such that $p\in y_1, y_2$ as follows: considering $y$, $y_1$ and $y_2$ as functions from $\omega$ to $\omega$, let $y_1$ and $y_2$ have $p$ as initial segment and be such that $y\backslash p = (y_0 \backslash p)\oplus (y_1\backslash p)$, where $\oplus$ is the operation of alternating digits between the reals.
		
		Then, in \(V[x, y_1, y_2]=V[x,y]\), we have 
		\[ B^I = \sum_{i=0}^{m-1} \tau_i^{1} S_i = \sum_{i=0}^{m-1} \tau_i^{2} S_i,\]
		where $\tau_i^{1}$ and $\tau_i^{2}$ are the interpretations of the name $\tau_i$ by $y_1$ and $y_2$ respectively, for each $i<m$. 
		In particular, $\tau_i^{1} \in V[y_1]$ and $\tau_i^{2} \in V[y_2]$.
		But on the other hand, by uniqueness of the coefficients \( \{c_i\}_{i=o}^{m-1}\), and taking into account that \(V[x,y_1], V[x,y_2] \subseteq V[x,y]\), we have that \(\tau_i^1= \tau_i^2 = \tau_i^{y}=c_i\) for $i=0, \dots, m-1$. 
		In particular, \(c_i \in \mathbb{R}^{V[y_1]} \cap \mathbb{R}^{V[y_2]} = \mathbb{R}^{V}. \)
		In other words, \(B^I = \sum_{i=0}^{m-1} c_i S_i\), where $c_i \in \mathbb{R}^V$. But this means \(0= B^I - \sum_{i=0}^{m-1} c_i S_i\), where $c_i \not = 0$ and the right hand side is not trivial (we chose $I$ such that $I\not \in \mathcal{J}$, i.e., $S_i\not = B^I$ for all $i=0,\dots, m-1$). This contradicts $B$ being algebraically independent over $\mathbb{R}^V$. 
	\end{Dem}
	
	\begin{Remark}
		Note that the same proof shows that if $B\subseteq \mathbb{R}^{V[x]}$ is \emph{linearly} independent over $\mathbb{R}^V$, then $B$ is also a linearly independent set over $\mathbb{R}^{V[y]}$. 
	\end{Remark}

	\bibliographystyle{siam}
	\bibliography{../../bibliography.bib}
	
\end{document}